\pdfoutput=1
\documentclass[11pt]{article}
\usepackage[T1]{fontenc}
\usepackage{lmodern}
\usepackage[margin=2.7cm]{geometry}
\usepackage{amsmath,amssymb,amsthm,mathtools}
\usepackage{microtype}
\usepackage{enumitem}
\usepackage{booktabs}
\usepackage{graphicx}
\usepackage{tikz}
\usetikzlibrary{arrows.meta}
\usepackage[round,authoryear]{natbib}
\usepackage[colorlinks=true,citecolor=blue,linkcolor=blue,urlcolor=blue]{hyperref}

\newcommand{\R}{\mathbb{R}}
\newcommand{\relint}{\operatorname{relint}}

\DeclareMathOperator{\diag}{diag}
\newcommand{\zenodolink}{\href{https://doi.org/10.5281/zenodo.23713443}{\texttt{doi:10.5281/zenodo.23713443}}}

\theoremstyle{plain}
\newtheorem{theorem}{Theorem}[section]
\newtheorem{proposition}[theorem]{Proposition}
\newtheorem{lemma}[theorem]{Lemma}
\newtheorem{corollary}[theorem]{Corollary}

\theoremstyle{definition}

\theoremstyle{remark}

\title{A five-species counterexample to the Hofbauer--Kon--Saito conjecture
on qualitative permanence}
\author{Miruna-Mihaela Beldiman\\[2pt]
\small Institute of Mathematical Statistics and\\
\small Applied Mathematics of the Romanian Academy\\
\small \texttt{miruna.m@gmail.com}
\and
Andrei-Dan Halanay\\[2pt]
\small Faculty of Mathematics and Computer Science,\\
\small University of Bucharest\\
\small \texttt{andrei.halanay@unibuc.ro}}
\date{October 2026}

\begin{document}
\maketitle

\begin{abstract}
Hofbauer, Kon and Saito showed that a Lotka--Volterra interaction matrix $A$ with
negative diagonal can be qualitatively permanent only if (C1) every cycle of
length $\ge2$ of $A$ is nonpositive and (C2) every negative cycle contains
exactly one negative entry. They conjectured that these conditions are also
sufficient, and proved this for $n\le3$. We give a $5\times5$ integer matrix
with seven nonzero off-diagonal entries that satisfies (C1) and (C2) but is
not qualitatively permanent. For an explicit integer growth vector, in which
one species is an obligate consumer as the definition allows, the system has
an interior equilibrium, and an asymptotically stable heteroclinic cycle on the boundary
of the orthant attracts an open set of interior orbits.
Adding non-interacting species gives counterexamples for every $n\ge5$. The
interactions form a one-way feedback loop together with a predator of two
prey. The verification is
elementary except for one heteroclinic connection, which lies in a face that
is not Volterra--Lyapunov stable. For it we prove a connection lemma based on
invasion graphs, which reduces the question to a smaller, stable face. The failure is robust. It persists under small
perturbations, including weak additional interactions, so every sign pattern
that satisfies (C1) and (C2) and contains ours is a counterexample, and it
occurs for an open cone of growth vectors. For $n\ge6$ there are also
counterexamples in which every species grows when alone. Finally, we announce
that the conjecture holds for $n=4$ and that, for $n=5$, four minimal sign
patterns characterise its failure.
\end{abstract}



\section{Introduction}

Consider the Lotka--Volterra system
\begin{equation}\label{eq:lv}
  \dot x_i=x_i\Bigl(r_i+\sum_{j=1}^na_{ij}x_j\Bigr),\qquad i=1,\dots,n,
\end{equation}
on $X=\R^n_{\ge0}$. It is \emph{permanent} if there is $\delta>0$ with
$\delta\le\liminf_{t\to\infty}x_i(t)\le\limsup_{t\to\infty}x_i(t)\le1/\delta$
for all $i$ and every solution with $x(0)\in\operatorname{int}X$. Let $Q_A$
be the set of matrices with the same sign pattern as $A$. Following
\citet[Def.~2.1]{hofbauer2008}, $A$ is \emph{qualitatively permanent} if
$\dot x_i=x_i\sum_j\tilde a_{ij}(x_j-x_j^*)$ is permanent for every
$\tilde A\in Q_A$ and every $x^*>0$. Equivalently, \eqref{eq:lv} is
permanent for every $\tilde A\in Q_A$ and every $r\in\R^n$ with
$r+\tilde Ax^*=0$ for some $x^*>0$, that is, for every $r$ for which the
system has an interior equilibrium. No sign condition is imposed on $r$.
A negative $r_i$ describes an obligate consumer, which declines when alone,
such as the predator of the classical Lotka--Volterra predator--prey model;
what decides whether all species can coexist at equilibrium is $x^*>0$, not
the sign of $r$. Much of the permanence literature concerns systems in which
every species can grow alone, $r\gg0$. That is a smaller class of systems;
see Remark~(v) below and Proposition~\ref{prop:rpos}.

A cycle of $A$ is a product $a_{i_1i_2}a_{i_2i_3}\cdots a_{i_ki_1}$ over
distinct indices, $k\ge2$.

\begin{theorem}[\citealp{hofbauer2008}, Thms.~2.3 and~2.7]
Let $a_{ii}<0$ for all $i$. If $A$ is qualitatively permanent, then
\begin{itemize}[leftmargin=3.2em]
\item[(C1)] every cycle of $A$ of length $\ge2$ is $\le0$, and
\item[(C2)] every negative cycle of $A$ contains exactly one negative entry.
\end{itemize}
For $n\le3$, (C1) and (C2) are also sufficient.
\end{theorem}

\citet[Conj.~2.4]{hofbauer2008} conjectured that (C1) and (C2) characterise
qualitative permanence for every $n$. They proved it for patterns of a
generalised food-chain form \citep[Thm.~2.6]{hofbauer2008}. To our knowledge
the conjecture has remained open.

\begin{theorem}\label{thm:main}
Let
\[
  A=\begin{pmatrix}
  -1&1&2&0&0\\ -2&-2&0&12&0\\ -1&0&-1&0&0\\ 0&0&0&-1&-24\\ 0&12&0&0&-1
  \end{pmatrix},\qquad r=(1,3,3,1,-2).
\]
\begin{enumerate}[label=(\alph*)]
\item $a_{ii}<0$ for all $i$, and $A$ satisfies (C1) and (C2).
\item \eqref{eq:lv} has the interior equilibrium
  $x^*\approx(2.390,\,0.170,\,0.610,\,0.177,\,0.034)$.
\item The boundary of $X$ contains a heteroclinic cycle
  \[
    \Gamma:\quad F_{13}\to F_{134}\to F_{124}\to F_{125}\to F_{13}
  \]
  through the equilibria $F_{13}=(\tfrac73,0,\tfrac23,0,0)$,
  $F_{134}=(\tfrac73,0,\tfrac23,1,0)$, $F_{124}=(\tfrac{17}4,\tfrac{13}4,0,1,0)$
  and $F_{125}=(\tfrac54,\tfrac14,0,0,1)$.
\item $\Gamma$ is asymptotically stable. In particular \eqref{eq:lv} is not
  permanent, and $A$ is not qualitatively permanent.
\end{enumerate}
Hence Conjecture~2.4 of \citet{hofbauer2008} is false for $n=5$.
\end{theorem}

Precisely, the sign pattern of $A$ satisfies (C1) and (C2), and for the
matrix $A$ in its class and the admissible $r$ above the system is not
permanent. So (C1) and (C2) are not sufficient for qualitative permanence.
They remain necessary \citep[Thm.~2.3]{hofbauer2008}.

\begin{corollary}\label{cor:alln}
The conjecture is false for every $n\ge5$.
\end{corollary}

\begin{proof}
For $n>5$ extend $A$ by $a_{ii}=-1$, $r_i=1$ and $a_{ij}=a_{ji}=0$ for
$i>5$. No new nonzero cycles appear, so (C1) and (C2) still hold, and
$(x^*,1,\dots,1)$ is an interior equilibrium. The new coordinates decouple
and converge to $1$ for every positive initial value. So
$\Gamma\times\{(1,\dots,1)\}$ is asymptotically stable, and it lies in the
boundary of $\R^n_{\ge0}$.
\end{proof}

\paragraph{The ecology.} We write $f_i(x)=r_i+(Ax)_i$ for the per-capita
growth rates.
\begin{itemize}
\item Species~1 eats species 2 and~3 ($a_{12},a_{13}>0$ and
  $a_{21},a_{31}<0$), and it also grows on its own ($r_1>0$).
\item Species~2 benefits from~4 ($a_{24}>0$), 5 benefits from~2 ($a_{52}>0$),
  and 5 harms~4 ($a_{45}<0$). None of these three effects is reciprocated:
  $a_{42}=a_{25}=a_{54}=0$.
\item Only species~5 declines when alone ($r_5<0$): it consumes~2.
\end{itemize}
One-sided effects such as commensalism ($a_{ij}>0=a_{ji}$) and amensalism
($a_{ij}<0=a_{ji}$) are common in ecological communities. Here they are not
essential, however. Weak reciprocal effects may be added as long as (C1) and
(C2) still hold (Corollary~\ref{cor:pattern}), and in the motif $M_4$ of
\S\ref{sec:disc} the harm is reciprocated by predation.
The three one-way links form a negative feedback loop
$2\to5\dashv4\to2$. It is the only cycle of $A$ of length $\ge3$, and it has
exactly one negative entry, as (C2) requires.

\begin{figure}[t]
\centering
\includegraphics[width=\textwidth]{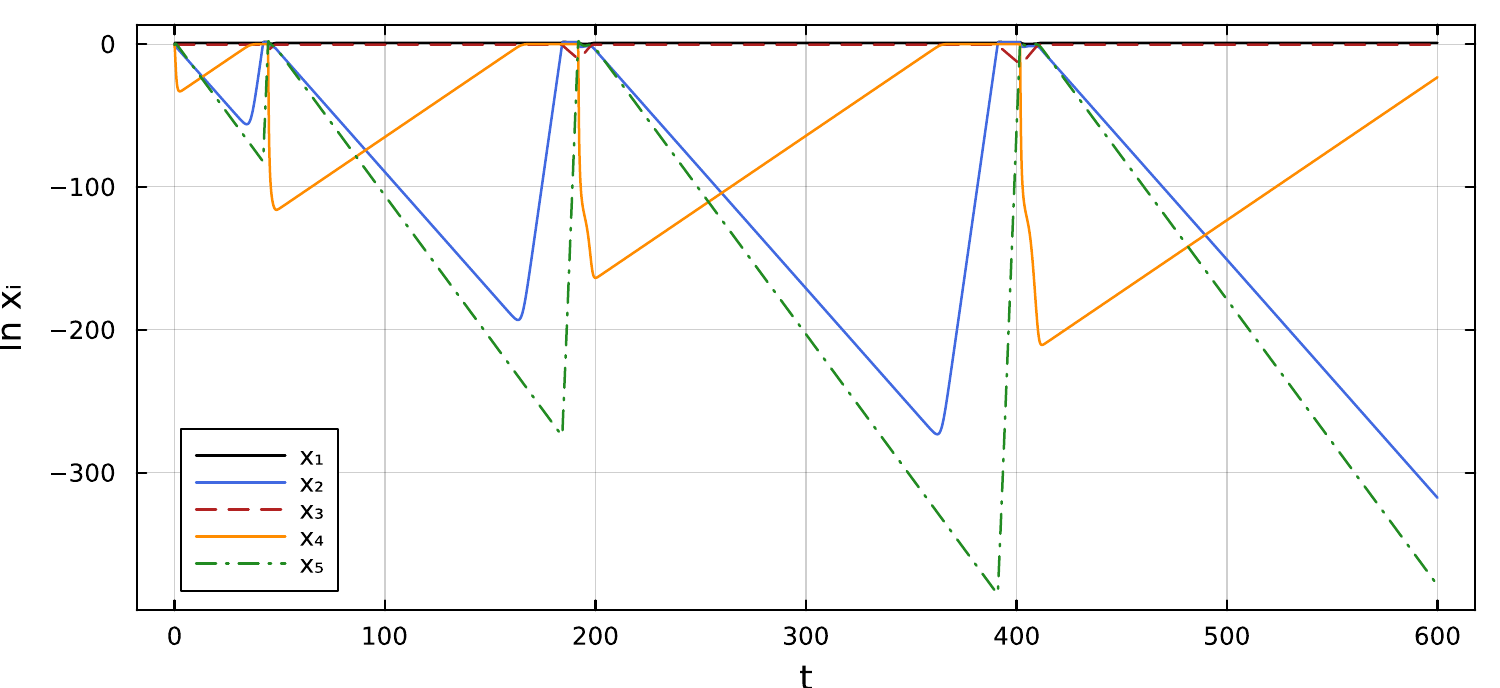}
\caption{$\ln x_i(t)$ along the interior orbit of \eqref{eq:lv} with
$x(0)=(2,\tfrac12,\tfrac12,\tfrac12,\tfrac12)$, for $A$ and $r$ of
Theorem~\ref{thm:main}. Species~1 persists. The others are repeatedly lost
and re-invade, and each excursion towards the boundary is deeper than the
previous one: the orbit approaches the heteroclinic cycle $\Gamma$.}
\label{fig:orbit}
\end{figure}

Along $\Gamma$, species~1 is always present:
\begin{itemize}
\item species~4 invades the predator--prey pair $\{1,3\}$;
\item species~2, helped by~4, invades and displaces~3;
\item species~5, helped by~2, invades and drives~4 out;
\item without~4, prey~3 returns and 2 and~5 are lost.
\end{itemize}
Figure~\ref{fig:orbit} shows an interior orbit following $\Gamma$ with
ever longer and deeper excursions.

\paragraph{Method.} Asymptotic stability of a simple heteroclinic cycle
through equilibria is decided by the characteristic matrix of
\citet{hofbauer1994}. That part is a finite computation with small
rationals, given in \S\ref{sec:proof}.

The connections need more care. Three of the four lie in faces where the
interaction matrix is Volterra--Lyapunov (VL) stable, and a Lyapunov function
gives them. The fourth lies in the face $\{x_3=0\}$, which contains the
feedback loop. The product $|a_{24}a_{45}a_{52}|=3456$ is $1728$ times the
product of the corresponding diagonal entries. The principal submatrix on $\{2,4,5\}$
is a one-way cyclic matrix with negative feedback, and such a matrix is
diagonally stable exactly when this ratio is below $\sec^3(\pi/3)=8$
\citep{arcak2006}. VL-stability passes to principal submatrices, so this
face is not VL-stable. For it we prove
Lemma~\ref{lem:morse}. In a face whose invasion graph is acyclic, every
interior orbit converges to the saturated equilibrium of the face, provided
the face of that equilibrium is VL-stable. So the invasion graph reduces the
question to a smaller face, which here does not contain the loop; it does
not replace Lyapunov functions there. The proof combines the Morse
decomposition of \citet{hofbauer2022} with time averages.

\paragraph{Remarks on the scope.}
\begin{enumerate}[label=(\roman*)]
\item The counterexample is robust: small perturbations of $A$ and $r$ keep
  an interior equilibrium and an asymptotically stable simple heteroclinic
  cycle with the same supports, because all hypotheses used below are open
  (Theorem~\ref{thm:robust}). Hence every pattern that satisfies (C1) and
  (C2) and contains the pattern of $A$ is a counterexample
  (Corollary~\ref{cor:pattern}). Section~\ref{sec:rset} also describes the
  set of growth vectors $r$ for which the proof applies, explains why the
  loop must be strong, and gives a normalised example with all entries of
  modulus at most~$6$ (Proposition~\ref{prop:normal}).
\item In the example $r_5<0$: species~5 is an obligate consumer of~2.
  Definition~2.1 of \citet{hofbauer2008} allows this, and the cycle of
  Theorem~\ref{thm:main} requires it, since $\lambda_5(F_{13})=r_5$ must be
  negative. For $n\ge6$ there are counterexamples with $r\gg0$
  (Proposition~\ref{prop:rpos}), by the same mechanism. For $n=5$ there is
  none: with $r\gg0$ every five-species system satisfying (C1) and (C2)
  with an interior equilibrium is permanent (Theorem~\ref{thm:five}(b),
  proved in a companion paper).
\item The pattern of $A$ does not have the form (2.2) of
  \citet[Thm.~2.6]{hofbauer2008}, under any relabelling. This is consistent
  with their theorem.
\item For $n=4$ the conjecture holds, for every $r$. We prove this in a
  companion paper by a computer-assisted argument (see \S\ref{sec:disc}).
  So $n=5$ is the smallest dimension in which it fails.
\item For $n=5$, at least seven nonzero off-diagonal entries are needed
  (companion paper). Every $5\times5$ pattern for which we found an
  attracting heteroclinic cycle contains one of four motifs, and the pattern
  of $A$ is one of them. We certify an attracting heteroclinic cycle for
  each motif (Proposition~\ref{prop:motifs}). In the companion paper we
  prove that the four motifs characterise the failure at $n=5$: a
  $5\times5$ sign pattern is qualitatively permanent if and only if it
  satisfies (C1) and (C2) and contains none of them
  (Theorem~\ref{thm:five}).
\end{enumerate}

\section{Preliminaries}

\paragraph{Faces and equilibria.} For $I\subseteq[n]=\{1,\dots,n\}$ let
$X_I=\{x\in X:x_j=0\text{ for }j\notin I\}$. It is invariant, and
$\relint X_I=\{x\in X_I:x_i>0\text{ for }i\in I\}$. If $A_{II}$ is
nonsingular and $\hat x_I=-A_{II}^{-1}r_I$ is positive, then $F_I$ denotes
the equilibrium with these coordinates on $I$ and $0$ elsewhere; it is
\emph{feasible}. For $j\notin I$ the \emph{invasion rate} of $j$ at $F_I$ is
$\lambda_j(F_I)=f_j(F_I)=r_j+(AF_I)_j$.

With the coordinates ordered as $(I,[n]\setminus I)$, the Jacobian at $F_I$ is
\[
  J(F_I)=\begin{pmatrix}\diag(\hat x_I)A_{II}&*\\0&\diag(\lambda_j(F_I))_{j\notin I}\end{pmatrix}.
\]
So the eigenvalues of $\diag(\hat x_I)A_{II}$ are \emph{internal} and the
invasion rates are \emph{external}. $F_I$ is \emph{internally stable} (IS)
if $\diag(\hat x_I)A_{II}$ is Hurwitz, and hyperbolic if moreover all its
invasion rates are nonzero.

For $I\subseteq K$, $F_I$ is a \emph{saturated equilibrium of the face}
$X_K$ if $\lambda_j(F_I)\le0$ for all $j\in K\setminus I$.

\paragraph{VL-stability.} A square matrix $M$ is \emph{VL-stable} if
$DM+M^\top D$ is negative definite for some positive diagonal matrix $D$.

\begin{lemma}\label{lem:VL}
Let $A_{KK}$ be VL-stable, and let $F_J$ ($J\subseteq K$) be a saturated
equilibrium of $X_K$. Then every orbit with $x(0)\in\relint X_K$ converges to
$F_J$. In particular $\{F_J\}$ is the only nonempty compact invariant subset
of $\relint X_J$.
\end{lemma}

\begin{proof}
This is the classical argument \citep[\S15.3]{hofbauer1998}. Let $\hat x=F_J$
and $V(x)=\sum_{i\in K}d_i(x_i-\hat x_i\log x_i)$ on $\relint X_K$. Then
\[
  \dot V=\sum_{i\in K\setminus J}d_ix_i\lambda_i(\hat x)
  +\tfrac12(x-\hat x)_K^\top(DA_{KK}+A_{KK}^\top D)(x-\hat x)_K
  \le-c|x-\hat x|^2
\]
for some $c>0$. The first sum is $\le0$ by saturation. So orbits are
bounded and $\int_0^\infty|x-\hat x|^2dt<\infty$, which gives
$x(t)\to\hat x$.

For the last claim apply this with $K=J$. On a compact invariant
$C\subseteq\relint X_J$, $V$ attains its maximum at some $y$. Then $V$ is
constant along the orbit of $y$ in $C$, so $y=\hat x$, and
$C=\{\hat x\}$.
\end{proof}

\paragraph{Heteroclinic cycles.} A \emph{boundary cycle of equilibria}
$\Gamma\subset\partial X$ consists of equilibria $F_1,\dots,F_m$ with supports
$I_1,\dots,I_m$ (indices mod $m$), together with orbits $\gamma_k$ with
$\alpha(\gamma_k)=F_k$ and $\omega(\gamma_k)=F_{k+1}$.

Its \emph{characteristic matrix} \citep{hofbauer1994} is the $m\times n$
matrix $C_{kj}=\lambda_j(F_k)$, with $C_{kj}=0$ for $j\in I_k$. Columns that
vanish identically belong to species present at every $F_k$ and are deleted.
$\Gamma$ is \emph{simple} if every row and every remaining column of $C$ has
exactly one positive entry. The positive entry of row $k$ is the
\emph{invader} $j_k$ at $F_k$. $\Gamma$ is \emph{hyperbolic}, resp.\ IS, if
every $F_k$ is.

Theorem~\ref{thm:hofbauer} below does not itself need simplicity. We use it
twice. First, at each $F_k$ exactly one absent species can invade, so the
unstable manifold of $F_k$ is a curve and $\Gamma$ is stable within
$\partial X$ (Proposition~\ref{prop:boundary}). Second, after ordering the
columns by the invaders, $C$ is square with its positive entries on the
diagonal, i.e.\ a Z-matrix. The condition of Theorem~\ref{thm:hofbauer}
then says that $C$ is not an M-matrix \citep[Cor.~1]{hofbauer1994}.

\begin{theorem}[\citealp{hofbauer1994}, Thm.~1(b); \citealp{hofbauer1998}, Thm.~17.5.1(b)]\label{thm:hofbauer}
Let $\Gamma$ be asymptotically stable within $\partial X$. If $Cp\gg0$ for
some $p\ll0$, then $\Gamma$ is asymptotically stable in $X$.
\end{theorem}

\begin{proposition}\label{prop:boundary}
A simple, hyperbolic, IS boundary cycle of equilibria is asymptotically
stable within $\partial X$.
\end{proposition}

This is stated parenthetically in \citet[Cor.~1]{hofbauer1994}. For
completeness, a proof is given in Appendix~\ref{app:boundary}.

\section{Two connection lemmas}

\begin{lemma}\label{lem:conn}
Let $F_I$ and $F_J$ be feasible, $J\not\subseteq I$, $K=I\cup J\ne[n]$, and
$|J\setminus I|=1$, say $J\setminus I=\{j\}$. Assume that $F_I$ is IS, that
$\lambda_j(F_I)>0$, that all other invasion rates at $F_I$ are negative,
and that every orbit of $\relint X_K$ converges to $F_J$. Then the unstable
manifold of $F_I$ is a curve, and one branch of it is an orbit
$\gamma\subset\relint X_K\subset\partial X$ with $\alpha(\gamma)=F_I$ and
$\omega(\gamma)=F_J$.
\end{lemma}

\begin{proof}
$J(F_I)$ has exactly one eigenvalue with positive real part, namely
$\lambda_j(F_I)$. Its eigenvector $v$ has $v_j\ne0$: otherwise $v$ would be
an eigenvector of an internal eigenvalue or of another invasion rate, all of
which have negative real part. Restricted to the invariant face $X_K$, the
unstable manifold $W^u(F_I)$ is a curve tangent to $v$, and one of its
branches has $x_j>0$. Near $F_I$ its coordinates in $I$ are positive too, so
this branch lies in $\relint X_K$. By hypothesis it converges to $F_J$.
\end{proof}

The hypothesis ``every orbit of $\relint X_K$ converges to $F_J$'' follows
from Lemma~\ref{lem:VL} when $A_{KK}$ is VL-stable and $F_J$ is saturated in
$X_K$. The next lemma is a dimensional reduction. It needs no Lyapunov
function on $X_K$, only acyclicity of the invasion graph of the face, and
it moves the stability requirement to the smaller face $X_J$ of the target,
where it is checked by Lemma~\ref{lem:VL}. In our application $A_{KK}$ is
not VL-stable, since its loop ratio exceeds the secant bound~$8$
\citep{arcak2006}, whereas $X_J$ does not contain the loop and is
qualitatively VL-stable. Its new ingredient is an inward transfer: the Morse decomposition that
an acyclic invasion graph gives on the boundary of the face
\citep{hofbauer2022}, combined with time averages and the uniqueness of
saturated equilibria, determines the limit of every interior orbit of the
face. Let $M$ be a square matrix. The \emph{invasion graph} of the
system $\dot x_i=x_i(r_i+(Mx)_i)$, $i\in K$ \citep{hofbauer2022}, has as
vertices the supports $L\subsetneq K$ of its feasible equilibria, and an
edge $S\to T$ whenever $S\neq T$, $\lambda_i(F_S)>0$ for all $i\in T\setminus
S$, and $\lambda_i(F_T)<0$ for all $i\in S\setminus T$.

\begin{lemma}[Connection lemma via invasion graphs]\label{lem:morse}
Let $K\subseteq[n]$ with $a_{ii}<0$ for $i\in K$, let $A_{KK}$ satisfy (C1),
and assume:
\begin{enumerate}[label=(\roman*)]
\item every invasion rate $\lambda_i(F_L)$, $i\in K\setminus L$, at every
  feasible $F_L$ with $L\subsetneq K$ is nonzero;
\item the invasion graph of the face system $(r_K,A_{KK})$ is acyclic;
\item $F_J$, $J\subsetneq K$, is feasible with $\lambda_i(F_J)<0$ for all
  $i\in K\setminus J$;
\item $\{F_J\}$ is the only nonempty compact invariant subset of
  $\relint X_J$. By Lemma~\ref{lem:VL}, this holds if $A_{JJ}$ is
  VL-stable. It holds for all magnitudes if $A_{JJ}$ has negative diagonal,
  satisfies (C1) and has no nonzero cycle of length $\ge3$. These sign
  conditions characterise qualitative VL-stability, that is, VL-stability of
  every matrix with the sign pattern of $A_{JJ}$
  \citep[(1.2)]{hofbauer2008}.
\end{enumerate}
Then every orbit with $x(0)\in\relint X_K$ converges to $F_J$.
\end{lemma}

The lemma uses (C1) and not (C2). (C1) enters through two facts: $-A_{KK}$
and all its principal submatrices are P-matrices
\citep[Lemma~6.2]{hofbauer2008}, and the face system has a compact
absorbing set \citep[Lemma~6.3]{hofbauer2008}. The first is a one-line
check. In the expansion of a principal minor $\det(-A_{LL})$, a permutation
contributes the product over its cycles $(i_1\dots i_k)$, $k\ge2$, of
$-a_{i_1i_2}\cdots a_{i_ki_1}\ge0$, and over its fixed points of
$-a_{ii}>0$. So every term is nonnegative, the identity contributes a
positive term, and negative 3-cycles, which (C2) allows, contribute
positively.

\begin{figure}[t]
\centering
\resizebox{\textwidth}{!}{%
\begin{tikzpicture}[>={Stealth[length=2mm]}, every node/.style={font=\small}]
\fill[gray!8] (0,0) rectangle (5,3);
\draw[->, thick, gray] (0.12,0) -- (2.38,0);
\draw[->, thick, gray] (2.62,0) -- (5,0) -- (5,1.38);
\draw[->, thick, gray] (0,0.12) -- (0,1.38);
\draw[->, thick, gray] (0,1.62) -- (0,3) -- (5,3) -- (5,1.62);
\fill (0,0) circle (2pt) node[below left] {$M_1$};
\fill (2.5,0) circle (2pt) node[below] {$M_2$};
\fill (0,1.5) circle (2pt) node[left] {$M_3$};
\fill (5,1.5) circle (2pt) node[right] {$M_q=\{F_J\}$};
\node at (2.5,-0.8) {$\partial X_K$ with its Morse sets};
\node at (1.0,2.45) {$\operatorname{relint}X_K$};
\fill (1.2,1.4) circle (1.3pt) node[below] {$z$};
\draw[->] (1.2,1.4) .. controls (2.0,2.4) and (3.8,2.3) .. (4.85,1.56);
\node[align=left, anchor=west, font=\footnotesize] at (6.9,1.5)
  {Step 2: $\omega(z)\subseteq\Gamma_0\subseteq\partial X_K$\\[2pt]
   Step 3: $\omega(z)$ lies in one Morse set $M_L$\\[2pt]
   Step 4: time averages give $L=J$\\[2pt]
   Step 5: (iv) gives $\omega(z)=\{F_J\}$};
\end{tikzpicture}}
\caption{Schematic of the proof of Lemma~\ref{lem:morse}. The acyclic
invasion graph orders the boundary Morse sets $M_1,\dots,M_q$; boundary
orbits (grey) run from lower to higher ones. An interior orbit cannot stay
away from the boundary, its limit set must lie in a single Morse set, and
time averages identify that set as the one of the saturated equilibrium
$F_J$.}
\label{fig:morse}
\end{figure}

\begin{proof}
The proof has five steps (Figure~\ref{fig:morse}): the setting of
\citet{hofbauer2022}; $\omega(z)$ lies in the boundary attractor; it lies
in a single Morse set; time averages identify that set; and (iv)
concludes. The face system on $X_K$ is again of the form \eqref{eq:lv}. By
\citet[Lemma~6.2]{hofbauer2008}, $-A_{KK}$ is a P-matrix, and so is every
principal submatrix. Fix $z\in\relint X_K$.
\begin{enumerate}[label=\arabic*.]
\item \emph{Assumptions of \citet{hofbauer2022}.} The face system has a
  compact absorbing set \citep[Lemma~6.3]{hofbauer2008}, so A1--A2 of
  \citet{hofbauer2022} hold.
  \begin{itemize}
  \item Let $\nu$ be an ergodic invariant measure with support $L$ (the set
    of species positive $\nu$-a.e.). Then $\int f_i\,d\nu=0$ for $i\in L$
    \citep[Lemma~1]{hofbauer2022}.
  \item Since each $f_i$ is affine, the mean $\bar x$ of $\nu$ solves
    $r_L+A_{LL}\bar x_L=0$, with $\bar x_L>0$. Hence $\bar x=F_L$ is feasible,
    and $\int f_i\,d\nu=\lambda_i(F_L)$ for $i\notin L$.
  \item With (i), this gives their assumption A3.
  \end{itemize}
\item \emph{$\omega(z)\subseteq\Gamma_0$, the part of the global attractor
  of the face system in $\partial X_K$.} A feasible equilibrium with
  support $K$ would be a saturated equilibrium of the face system different
  from $F_J$. This contradicts uniqueness of saturated equilibria for
  P-matrices \citep[Lemma~6.4]{hofbauer2008}, \citep[Thm.~15.4.5]{hofbauer1998}.
  So the face system has no equilibrium in $\relint X_K$. Then
  $\relint X_K$ contains no $\alpha$- or $\omega$-limit point of any orbit
  \citep[Thm.~5.2.1]{hofbauer1998}. This excludes periodic orbits and any
  other attractor in $\relint X_K$, and gives $\omega(z)\subseteq\partial X_K$.
  Since $\omega(z)$ is nonempty, compact and invariant, it lies in the
  global attractor, hence in $\Gamma_0$.
\item \emph{One Morse set.} For each vertex $L$ of the
  invasion graph let $M_L$ be the maximal compact invariant subset of
  $\relint X_L$. By (ii) and \citet[Lemma~4]{hofbauer2022}, these sets form a
  Morse decomposition of $\Gamma_0$. Order them as $M_1,\dots,M_q$ so that
  every other orbit of $\Gamma_0$ has $\alpha$-limit in some $M_i$ and
  $\omega$-limit in some $M_k$ with $k>i$. The union $A_j$ of
  $M_j,\dots,M_q$ and of the orbits of $\Gamma_0$ with both limit sets in
  them is an attractor of the flow on $\Gamma_0$ \citep{conley1978}, and
  $A_1=\Gamma_0$.

  By Step~2, $\omega(z)\subseteq\Gamma_0$. It is compact, invariant and internally chain
  transitive: for all $\varepsilon,T>0$, any two of its points are joined by
  an $(\varepsilon,T)$-chain inside it. Such a set, if it meets the basin of an
  attractor, is contained in the attractor \citep[see][Thm.~3.1, stated for maps]{hirsch2001}. For
  completeness, here is the argument.
  \begin{itemize}
  \item Let $\mathcal A$ be an attractor of the flow $\varphi$ on $\Gamma_0$
    and $\Omega$ such a set, meeting its basin. Then $\Omega$ meets
    $\mathcal A$, since it contains the $\omega$-limits of its points.
  \item $\mathcal A$ has a trapping region \citep{conley1978}: a compact,
    forward invariant neighbourhood $N$ in $\Gamma_0$ and $T>0$ with
    $\varphi_T(N)\subseteq\operatorname{int}N$ and
    $\bigcap_{t\ge0}\varphi_t(N)=\mathcal A$.
  \item If $\Omega\subseteq N$, then $\Omega=\varphi_t(\Omega)\subseteq\varphi_t(N)$
    for all $t$, so $\Omega\subseteq\mathcal A$.
  \item Otherwise pick $a\in\Omega\cap\mathcal A$ and $y\in\Omega\setminus N$, and
    let $\varepsilon$ be the distance from $\varphi_T(N)$ to
    $\Gamma_0\setminus\operatorname{int}N$. Since $\varphi_t(N)\subseteq\varphi_T(N)$
    for $t\ge T$, every $(\varepsilon,T)$-chain that starts in $N$ stays in $N$,
    so no such chain leads from $a$ to $y$. This is a contradiction.
  \end{itemize}
  Let $j$ be maximal with
  $\omega(z)\cap A_j\ne\emptyset$. Then $\omega(z)\subseteq A_j$ and
  $\omega(z)\cap A_{j+1}=\emptyset$. A point $y\in\omega(z)\setminus M_j$
  would lie in $A_{j+1}$, or have $\omega(y)\subseteq A_{j+1}$ with
  $\omega(y)\subseteq\omega(z)$. Both are impossible. Hence
  $\omega(z)\subseteq M_j$. Let $L$ be its support.
\item \emph{$L=J$.} Let $\mu$ be a weak$^*$ limit of
  $\frac1t\int_0^t\delta_{x(s)}ds$ along some $t_k\to\infty$, where $x(\cdot)$
  is the orbit of $z$. Then $\mu$ is invariant and supported on
  $\omega(z)\subseteq\relint X_L$.
  \begin{itemize}
  \item For $i\in K$,
    $\int f_i\,d\mu=\lim_k\frac1{t_k}\log\frac{x_i(t_k)}{x_i(0)}\le0$,
    because $x_i$ is bounded.
  \item All ergodic components of $\mu$ have support $L$, so by Step~1,
    $\int f_i\,d\mu=\lambda_i(F_L)$ for $i\in K\setminus L$.
  \end{itemize}
  So $F_L$ is feasible with $\lambda_i(F_L)\le0$ on $K\setminus L$, i.e.\
  saturated in $X_K$. By uniqueness, $L=J$.
\item \emph{Conclusion.} $\omega(z)$ is a nonempty compact invariant subset
  of $\relint X_J$. By (iv), $\omega(z)=\{F_J\}$.\qedhere
\end{enumerate}
\end{proof}

\paragraph{The case $K=[n]$.} The lemma also applies to the whole orthant,
$K=[n]$, but its conclusion is then extinction, not permanence. By
uniqueness of saturated equilibria, hypothesis (iii) with $J\subsetneq K$
means that the system has no interior equilibrium, and every interior orbit
converges to the boundary equilibrium $F_J$. When an interior equilibrium
exists, which is the setting of qualitative permanence, it is the saturated
equilibrium, so $J=K$ and the lemma does not apply. An acyclic invasion
graph then gives permanence instead, by \citet[Thm.~1]{hofbauer2022}; this
is how the result for $n=4$ in \S\ref{sec:disc} is obtained. So the lemma
says nothing for or against qualitative permanence. It is a tool for
connections inside proper faces.

\paragraph{Relation to \citet{hofbauer2022}.} Steps~1 and~3 use their
results: the Morse decomposition of the boundary attractor given by an
acyclic invasion graph (their Lemma~4), and the invasion rates of ergodic
measures (their Lemma~1). They use these to decide permanence, i.e.\
whether the boundary repels. Lemma~\ref{lem:morse} draws a different
conclusion: every interior orbit of a face converges to a prescribed
boundary equilibrium, although the face is not VL-stable. Step~2 is a
standard consequence of \citet[Thm.~5.2.1]{hofbauer1998}. The new part is
the combination in the remaining steps:
\begin{itemize}
\item the face has no interior equilibrium (Step~2);
\item the Morse set that contains $\omega(z)$ is identified by time averages
  and the uniqueness of saturated equilibria for P-matrices (Step~4);
\item hypothesis (iv) is applied inside that Morse set (Step~5).
\end{itemize}
To our knowledge, this use of invasion graphs to prove heteroclinic
connections is new.

\section{Proof of Theorem~\ref{thm:main}}\label{sec:proof}

All numbers below are exact; they can be checked by hand. The script
\texttt{hks\_note\_checks.jl} redoes all of them in rational arithmetic.

\paragraph{(a) Conditions (C1) and (C2).} The nonzero off-diagonal entries of
$A$ are $a_{12}$, $a_{13}$, $a_{21}$, $a_{24}$, $a_{31}$, $a_{45}$ and $a_{52}$. Index~1 lies only
on the 2-cycles through $(1,2)$ and $(1,3)$. Among $\{2,3,4,5\}$ the only
edges are $2\to4\to5\to2$. Hence the nonzero cycles of $A$ are exactly
\[
  a_{12}a_{21}=-2,\qquad a_{13}a_{31}=-2,\qquad a_{24}a_{45}a_{52}=-3456,
\]
each negative with exactly one negative entry.

\paragraph{(b) The interior equilibrium.} The equations $r+Ax=0$ give
$x_3=3-x_1$, $x_5=12x_2-2$ and $x_4=1-24x_5$. Then $3x_1-x_2=7$ and
$2x_1+3458x_2=591$, so $x_2=1759/10376$ and
$x^*=(24797,1759,6331,1832,356)/10376$.

\paragraph{(c) Equilibria and invasion rates.} Solving $r_I+A_{II}x_I=0$
gives the equilibria of Theorem~\ref{thm:main}. Their invasion rates are as
follows; the invader is in bold.
\begin{center}
\begin{tabular}{lccccc}
\toprule
 & $\lambda_1$ & $\lambda_2$ & $\lambda_3$ & $\lambda_4$ & $\lambda_5$\\
\midrule
$F_{13}=(\frac73,0,\frac23,0,0)$ & & $-\frac53$ & & $\mathbf 1$ & $-2$\\[2pt]
$F_{134}=(\frac73,0,\frac23,1,0)$ & & $\mathbf{\frac{31}3}$ & & & $-2$\\[2pt]
$F_{124}=(\frac{17}4,\frac{13}4,0,1,0)$ & & & $-\frac54$ & & $\mathbf{37}$\\[2pt]
$F_{125}=(\frac54,\frac14,0,0,1)$ & & & $\mathbf{\frac74}$ & $-23$ & \\
\bottomrule
\end{tabular}
\end{center}

\emph{Internal stability.} In each case $\diag(\hat x_I)A_{II}$ is
block-triangular. Its blocks are:
\begin{itemize}
\item a $2\times2$ predator--prey block on $\{1,3\}$ or $\{1,2\}$, of the form
  $\bigl(\begin{smallmatrix}-x_1&x_1a_{1k}\\x_ka_{k1}&x_ka_{kk}\end{smallmatrix}\bigr)$
  with $a_{1k}a_{k1}<0$, so negative trace and positive determinant;
\item for species 4 or~5, a $1\times1$ block equal to $-1$.
\end{itemize}
Indeed, in $F_{134}$ and $F_{124}$ row~4 of $A_{II}$ vanishes off the diagonal,
and in $F_{125}$ column~5 does. So all four equilibria are IS. Since all
invasion rates are nonzero, they are also hyperbolic.

\emph{Simplicity.} The column of species~1 vanishes (species~1 is present
throughout). In the remaining columns $2,3,4,5$ and in each row there is
exactly one positive entry. So, once the connections below exist, $\Gamma$
is a simple, hyperbolic, IS boundary cycle with invaders $4,2,5,3$.

\paragraph{(d) The connections.} In each step $F_I\to F_J$, the set
$J\setminus I$ consists of the invader alone, and $K=I\cup J\ne[5]$. By
Lemma~\ref{lem:conn} it suffices to show that every orbit of $\relint X_K$
converges to $F_J$.

\emph{$F_{13}\to F_{134}$, $K=\{1,3,4\}$; $F_{134}\to F_{124}$,
$K=\{1,2,3,4\}$; $F_{125}\to F_{13}$, $K=\{1,2,3,5\}$.} Here $A_{KK}$ is
VL-stable, with
\[
  D=\diag(1,2,1),\qquad D=\diag(1,\tfrac12,2,10),\qquad D=\diag(1,\tfrac12,2,\tfrac1{48})
\]
respectively. The predator--prey entries cancel in $DA+A^\top D$, and
$-(DA_{KK}+A_{KK}^\top D)$ is
\[
  \diag(2,4,2),\qquad
  \begin{pmatrix}2&0&0&0\\0&2&0&-6\\0&0&4&0\\0&-6&0&20\end{pmatrix},\qquad
  \begin{pmatrix}2&0&0&0\\0&2&0&-\frac14\\0&0&4&0\\0&-\frac14&0&\frac1{24}\end{pmatrix},
\]
all positive definite ($2\cdot20>36$ and $2\cdot\frac1{24}>\frac1{16}$). The
target is saturated in $X_K$:
\begin{itemize}
\item $F_{134}$ has full support in $K=\{1,3,4\}$;
\item at $F_{124}$, $\lambda_3=-\frac54$;
\item at $F_{13}$, $\lambda_2=-\frac53$ and $\lambda_5=-2$.
\end{itemize}
Lemma~\ref{lem:VL} applies.

\emph{$F_{124}\to F_{125}$, $K=\{1,2,4,5\}$.} This face contains the loop,
and $A_{KK}$ is not VL-stable (\S1). We apply Lemma~\ref{lem:morse}; $A_{KK}$
satisfies (C1) by (a).
\begin{itemize}
\item[(i)] The feasible equilibria of the face system with support
  $L\subsetneq K$ and their invasion rates are as follows.
  \begin{center}\small
  \begin{tabular}{lll@{\qquad}lll}
  \toprule
  $L$ & $F_L$ on $L$ & rates & $L$ & $F_L$ on $L$ & rates\\
  \midrule
  $\emptyset$ & -- & $\lambda_{1,2,4,5}=1,3,1,-2$ &
  $\{1,4\}$ & $(1,1)$ & $\lambda_{2,5}=13,-2$\\
  $\{1\}$ & $(1)$ & $\lambda_{2,4,5}=1,1,-2$ &
  $\{2,4\}$ & $(\frac{15}2,1)$ & $\lambda_{1,5}=\frac{17}2,88$\\
  $\{2\}$ & $(\frac32)$ & $\lambda_{1,4,5}=\frac52,1,16$ &
  $\{2,5\}$ & $(\frac32,16)$ & $\lambda_{1,4}=\frac52,-383$\\
  $\{4\}$ & $(1)$ & $\lambda_{1,2,5}=1,15,-2$ &
  $\{1,2,4\}$ & $(\frac{17}4,\frac{13}4,1)$ & $\lambda_5=37$\\
  $\{1,2\}$ & $(\frac54,\frac14)$ & $\lambda_{4,5}=1,1$ &
  $\{1,2,5\}$ & $(\frac54,\frac14,1)$ & $\lambda_4=-23$\\
  \bottomrule
  \end{tabular}
  \end{center}
  All rates are nonzero. (The remaining supports, $\{5\}$, $\{1,5\}$,
  $\{4,5\}$, $\{1,4,5\}$ and $\{2,4,5\}$, give no feasible equilibrium.)
\item[(ii)] The invasion graph is acyclic (Figure~\ref{fig:graphs}(b)). The argument uses the sign
  pattern and four values of $x_2$, not the rates in the table. In
  $X_K$,
  \begin{gather*}
    f_1=1-x_1+x_2,\qquad f_2=3-2x_1-2x_2+12x_4,\\
    f_4=1-x_4-24x_5,\qquad f_5=-2+12x_2-x_5.
  \end{gather*}
  A species $i$ is lost along an edge $S\to T$ if $i\in S\setminus T$, and
  this needs $\lambda_i(F_T)<0$.
  \begin{itemize}
  \item \emph{Species 1 is never lost:} $\lambda_1(F_T)=1+x_2>0$.
  \item \emph{Species 2 is never lost.} If $2\notin T$, then $x_1\le1$ at
    $F_T$, because $f_1=1-x_1$ there. So $\lambda_2(F_T)\ge3-2=1>0$.
  \item \emph{Species 5 is never lost.} Every vertex containing 5 contains
    2, since $f_5<0$ when $x_2=0$. So if 5 is lost, either 2 is lost as well,
    which is excluded, or $2\in T$ and $5\notin T$. In the second case $T$
    is $\{2\}$, $\{1,2\}$, $\{2,4\}$ or $\{1,2,4\}$, where
    $x_2=\frac32,\frac14,\frac{15}2,\frac{13}4$. So
    $\lambda_5(F_T)=-2+12x_2\ge1>0$.
  \end{itemize}
  So $S\cap\{1,2,5\}\subseteq T\cap\{1,2,5\}$ along every edge. Along a
  cycle of the graph this set is therefore constant, say $P$, and consecutive
  vertices differ only in species~4. The cycle would then be
  $P\to P\cup\{4\}\to P$, which needs both $\lambda_4(F_P)>0$ and
  $\lambda_4(F_P)<0$.
\item[(iii)] $\lambda_4(F_{125})=-23<0$.
\item[(iv)] $A_{\{1,2,5\}}$ is VL-stable with $D=\diag(1,\frac12,\frac1{48})$:
  \[
    -(DA_{\{1,2,5\}}+A_{\{1,2,5\}}^\top D)=\begin{pmatrix}2&0&0\\0&2&-\frac14\\0&-\frac14&\frac1{24}\end{pmatrix}
  \]
  is positive definite.
\end{itemize}

\paragraph{(e) Stability.} With rows $F_{13},F_{134},F_{124},F_{125}$ and
columns ordered as the invaders $4,2,5,3$,
\[
  C=\begin{pmatrix}c_1&-\alpha&-\beta&0\\0&c_2&-\gamma&0\\0&0&c_3&-\delta\\-\varepsilon&0&0&c_4\end{pmatrix}
   =\begin{pmatrix}1&-\frac53&-2&0\\0&\frac{31}3&-2&0\\0&0&37&-\frac54\\-23&0&0&\frac74\end{pmatrix}.
\]
As for every simple cycle, $C$ is a Z-matrix: positive diagonal, other
entries $\le0$. Its leading principal minors $1$, $\frac{31}3$ and
$\frac{1147}3$ are positive, so $C$ is an M-matrix exactly when $\det C>0$
\citep[Cor.~1]{hofbauer1994}. The digraph of the negative entries has two
cycles, through the rows $1\to2\to3\to4\to1$ and $1\to3\to4\to1$. Both
pass through row~1, so
\[
  \det C=c_1c_2c_3c_4\,(1-g),\qquad
  g=\frac{\varepsilon}{c_4}\,\frac{\delta}{c_3}
    \Bigl(\frac{\beta}{c_1}+\frac{\alpha\gamma}{c_1c_2}\Bigr)
   =\frac{8280}{8029}>1,
\]
and $\det C=-\frac{251}{12}<0$. So $C$ is not an M-matrix, and there is
$p\ll0$ with $Cp\gg0$. Explicitly,
\[
  p=-(200,17,88,2620),\qquad Cp=\bigl(\tfrac{13}3,\tfrac13,19,15\bigr)\gg0,
\]
with $p_1=-1$ for the deleted column. By Proposition~\ref{prop:boundary},
$\Gamma$ is asymptotically stable within $\partial X$, and by
Theorem~\ref{thm:hofbauer} it is asymptotically stable in $X$.

The number $g$ is a loop gain. The ratios in it compare the decay rates of
absent species with the growth rates of the invaders, along the two cycles
of dependencies, and $\Gamma$ attracts because decay wins. Equivalently,
writing $C=D-N$ with $D$ the diagonal part, the spectral radius of
$D^{-1}N$ is $\approx1.0099>1$. For this $r$ the margin $g-1\approx0.03$ is
small. It is not small in general: \S\ref{sec:rset} gives growth vectors
with $g\approx4.6$, a certified box of them, and a normalised example with
$g\approx1.34$.

\begin{figure}[t]
\centering
\resizebox{\textwidth}{!}{%
\begin{tikzpicture}[>={Stealth[length=2mm]}, every node/.style={font=\small},
  v/.style={draw, rounded corners=2pt, inner sep=2pt, minimum width=7mm}]
\begin{scope}
\node (a) at (0,1.8) {$F_{13}$};
\node (b) at (3.4,1.8) {$F_{134}$};
\node (c) at (3.4,0) {$F_{124}$};
\node (d) at (0,0) {$F_{125}$};
\draw[->] (a) -- node[above] {$+4$} (b);
\draw[->] (b) -- node[right] {$+2,\,-3$} (c);
\draw[->] (c) -- node[below] {$+5,\,-4$} (d);
\draw[->] (d) -- node[left] {$+3,\,-2,\,-5$} (a);
\node at (1.7,-1.2) {(a)};
\end{scope}
\begin{scope}[xshift=6.3cm]
\node[v] (e) at (0,0) {$\emptyset$};   \node[v] (e4) at (0,1.8) {$4$};
\node[v] (o) at (1.4,0) {$1$};          \node[v] (o4) at (1.4,1.8) {$14$};
\node[v] (t) at (2.8,0) {$2$};          \node[v] (t4) at (2.8,1.8) {$24$};
\node[v] (ot) at (4.2,0) {$12$};        \node[v] (ot4) at (4.2,1.8) {$124$};
\node[v] (tf) at (5.6,0) {$25$};
\node[v] (otf) at (7.0,0) {$125$};
\draw[->] (e) -- (o); \draw[->] (e) to[bend right=28] (t); \draw[->] (e) -- (e4);
\draw[->] (o) to[bend right=48] (ot); \draw[->] (o) -- (o4);
\draw[->] (t) -- (ot); \draw[->] (t) -- (t4);
\draw[->] (ot) -- (ot4); \draw[->] (tf) -- (otf);
\draw[->] (e4) -- (o4); \draw[->] (e4) to[bend left=28] (t4);
\draw[->] (o4) to[bend left=48] (ot4); \draw[->] (t4) -- (ot4);
\draw[->] (t4) -- (tf); \draw[->] (ot4) -- (otf);
\node at (3.5,-1.2) {(b)};
\end{scope}
\end{tikzpicture}}
\caption{(a) The heteroclinic cycle $\Gamma$ of Theorem~\ref{thm:main}. Each
arrow is labelled with the invading species ($+$) and the species lost
($-$). (b) The invasion graph of the face system on $K=\{1,2,4,5\}$, with
supports written as strings. The columns are the values of
$S\cap\{1,2,5\}$, and the top row contains species~4. Shown are the 15 edges
that are not shortcuts of two-step paths; each of the other 11 edges joins
two vertices already joined by such a path. Every edge points to the right,
or straight up: species 1, 2 and~5 are never lost, as \S\ref{sec:proof}(d)(ii)
shows, so the graph is acyclic.}
\label{fig:graphs}
\end{figure}

An asymptotically stable $\Gamma\subset\partial X$ attracts every point of a
neighbourhood, including interior points, so \eqref{eq:lv} is not permanent.
Together with (b), this shows that $A$ is not qualitatively permanent.
\hfill$\square$

\section{Robustness and the set of growth rates}\label{sec:rset}

Qualitative permanence asks for permanence for every $\tilde A\in Q_A$ and
every admissible $r$, so one pair suffices to refute it. The failure is not
confined to a thin set of parameters, however. It persists under all small
perturbations, including new weak interactions, and it occurs for an open
cone of growth vectors.

For $n\ge5$ let $E=\{6,\dots,n\}$, let $A_n$ be $A$ extended by $a_{ii}=-1$
and $a_{ij}=a_{ji}=0$ for $i\in E$, and let $r_n=(r,1,\dots,1)$, as in
Corollary~\ref{cor:alln}.

\begin{theorem}\label{thm:robust}
Let $n\ge5$. There is $\varepsilon>0$ with the following property. Let
$\tilde A\in\R^{n\times n}$ satisfy (C1), and let $\tilde r\in\R^n$, with
$|\tilde a_{ij}-(A_n)_{ij}|<\varepsilon$ and $|\tilde r_i-(r_n)_i|<\varepsilon$
for all $i,j$. Then the system with $(\tilde A,\tilde r)$ has an interior
equilibrium and an asymptotically stable heteroclinic cycle in $\partial X$,
through equilibria with supports $\{1,3\}\cup E$, $\{1,3,4\}\cup E$,
$\{1,2,4\}\cup E$ and $\{1,2,5\}\cup E$. In particular it is not permanent.
\end{theorem}

\begin{proof}
First, the proof of \S\ref{sec:proof} applies to $(A_n,r_n)$ with every
support $I$ replaced by $I\cup E$.
\begin{itemize}
\item The species of $E$ are decoupled, with equilibrium value $1$. The
  rates of the species $1,\dots,5$ are unchanged, and the Jacobian blocks of
  $E$ are $-1$.
\item The VL faces and their certificates extend by $d_i=1$ on $E$.
\item In the face $K'=\{1,2,4,5\}\cup E$, a species of $E$ that is absent
  from a vertex has rate $1>0$, so it is never lost. Along a cycle of the
  invasion graph the set of species of $E$ is therefore constant, and
  (d)(ii) excludes the cycle. All rates are nonzero, by the table in (d)(i).
\item In the characteristic matrix, the columns of $1$ and of $E$ vanish.
\end{itemize}
Every hypothesis used there is open in $(\tilde A,\tilde r)$, as long as
$\tilde A$ satisfies (C1), which Lemma~\ref{lem:morse} requires.
\begin{itemize}
\item Feasibility, positivity of $x^*$ and the signs of the invasion rates
  are strict inequalities between continuous functions. The relevant
  principal submatrices stay nonsingular, since $-\tilde A$ is a P-matrix
  \citep[Lemma~6.2]{hofbauer2008}.
\item Internal stability (Hurwitz), VL-stability (with the same $D$) and
  $Cp\gg0$ (with the same $p$) are open.
\item Hypotheses (i) and (ii) of Lemma~\ref{lem:morse} on $K'$ persist,
  because the invasion graph does not change near $(A_n,r_n)$. A support
  $L\subsetneq K'$ with a negative coordinate of $-A_{LL}^{-1}r_L$ stays
  infeasible. If $-A_{LL}^{-1}r_L\ge0$ had a zero coordinate, it would be a
  feasible equilibrium with a smaller support and a zero invasion rate,
  which (i) excludes. So the feasible supports, and the signs of their rates,
  are locally constant.
\end{itemize}
\end{proof}

\begin{corollary}\label{cor:pattern}
Let $n\ge5$, and let $P$ be an $n\times n$ sign pattern with negative
diagonal that satisfies (C1) and (C2). Suppose that, after relabelling, $P$
has the sign of $a_{ij}$ at every nonzero entry $a_{ij}$ of $A$. Then $P$ is
not qualitatively permanent.
\end{corollary}

\begin{proof}
Let $\tilde A=A_n+\varepsilon S$, where $S_{ij}=\operatorname{sign}P_{ij}$
at the off-diagonal positions where $A_n$ is zero, and $S_{ij}=0$
elsewhere. For $\varepsilon>0$ small, $\tilde A$ has the sign pattern $P$,
so it satisfies (C1). Theorem~\ref{thm:robust} applies with $\tilde r=r_n$.
\end{proof}

So the counterexample is a sign motif, not a single matrix. Every pattern
that satisfies (C1) and (C2) and contains the pattern of $A$ fails,
whatever its other interactions. Explicitly, ten of the $13$ zero entries
of $A$ can be filled, each with one sign only: $a_{41},a_{42},a_{25}<0$ and
$a_{51},a_{32},a_{23},a_{43},a_{54},a_{15},a_{35}>0$. The entries $a_{53}$,
$a_{14}$ and $a_{34}$ cannot be filled at all. Altogether $223$ of the $3^{13}$
sign fillings keep (C1) and (C2), with up to six new entries, and
Corollary~\ref{cor:pattern} applies to each (script
\texttt{hks\_fill\_signs.jl}). In particular, each reciprocal entry of the
loop, $a_{42}$, $a_{25}$ or $a_{54}$, may be filled only with the
predator--prey sign, and at most two of them at once: all three would close
the positive cycle $a_{25}a_{54}a_{42}>0$. Since the entries of $A$ are bounded, Theorem~\ref{thm:robust}
holds equally for relative perturbations, each entry multiplied by a factor
close to~$1$.

\paragraph{Scale invariance and the strength of the loop.} Rescaling the
densities, $x_j\mapsto s_jx_j$, and the per-capita rates, $f_i\mapsto c_if_i$
with $c_i,s_j>0$, replaces $(r,A)$ by $(\diag(c)r,\diag(c)A\diag(s))$. This
preserves the sign pattern, feasibility, the signs of all invasion rates,
VL-stability, invasion graphs, and the criterion of
Theorem~\ref{thm:hofbauer} (the columns of $C$ are multiplied by positive
factors). It also preserves internal stability for the block structures of
\S\ref{sec:proof}. The magnitudes of $A$ therefore matter only through its
cycle ratios $|a_{i_1i_2}\cdots a_{i_ki_1}|/|a_{i_1i_1}\cdots a_{i_ki_k}|$,
which are $1$, $2$ and the loop ratio
$\ell=|a_{24}a_{45}a_{52}|/|a_{22}a_{44}a_{55}|=1728$. The value of $\ell$
cannot be small.

\begin{proposition}\label{prop:loop}
Let $\tilde A\in Q_A$ have loop ratio $\ell<8$. Then $\tilde A$ is
VL-stable, and \eqref{eq:lv} is permanent for every $r$ with an interior
equilibrium. So every counterexample in $Q_A$ has $\ell\ge8$.
\end{proposition}

\begin{proof}
Take $d_1=1$, and choose $d_2,d_3$ so that the predator--prey entries cancel
in $D\tilde A+\tilde A^\top D$. Then $-(D\tilde A+\tilde A^\top D)$ is block
diagonal, with positive blocks for species 1 and~3. The block on
$\{2,4,5\}$ is $-(D'\tilde A_{245}+\tilde A_{245}^\top D')$ with
$D'=\diag(d_2,d_4,d_5)$, and $\tilde A_{245}$ is a one-way cyclic matrix
with negative diagonal and negative feedback. The secant criterion
\citep{arcak2006} says that such a matrix is diagonally stable, i.e.\ some
$D'$ makes this block positive definite, exactly when
$\ell<\sec^3(\pi/3)=8$. Rescaling these weights by a common factor makes
$d_2$ the value already fixed, so $\tilde A$ is VL-stable. With an interior
equilibrium $x^*$, the function $V$ of Lemma~\ref{lem:VL} shows that every
interior orbit converges to $x^*$ \citep[\S15.3]{hofbauer1998}, so the
system is permanent.
\end{proof}

So the large entries of $A$ are not an accident of the search. After
normalising the diagonal to $-1$, the loop must have strength
$|a_{24}a_{45}a_{52}|\ge8$. In particular no counterexample in $Q_A$ has
all loop entries of modulus less than~$2$. With the cycle and the proof of
Theorem~\ref{thm:main}, a random search over $\tilde A$ and $r$ finds
counterexamples with $\ell$ down to about $14$; the predator--prey ratios
then grow (script \texttt{hks\_rcone.jl}).

\begin{proposition}\label{prop:normal}
Let
\[
  A'=\begin{pmatrix}
  -1&3&6&0&0\\ -3&-1&0&3&0\\ -6&0&-1&0&0\\ 0&0&0&-1&-3\\ 0&3&0&0&-1
  \end{pmatrix}\in Q_A,\qquad r'=(1,9,19,1,-1).
\]
Then the conclusions of Theorem~\ref{thm:main} hold for $(A',r')$, with the
same cycle, the interior equilibrium $x^*=(3283,432,157,292,251)/1045$,
\begin{gather*}
  F_{13}=\bigl(\tfrac{115}{37},0,\tfrac{13}{37},0,0\bigr),\qquad
  F_{134}=\bigl(\tfrac{115}{37},0,\tfrac{13}{37},1,0\bigr),\\
  F_{124}=\bigl(\tfrac{37}{10},\tfrac9{10},0,1,0\bigr),\qquad
  F_{125}=\bigl(\tfrac{14}5,\tfrac35,0,0,\tfrac45\bigr),
\end{gather*}
loop ratio $\ell=27$ and loop gain $g=\frac{8288}{6171}$.
\end{proposition}

\begin{proof}
$A'$ has the diagonal $-1$, symmetric predator--prey pairs and all loop
entries of modulus~$3$. The argument of \S\ref{sec:proof} applies verbatim:
the faces other than $\{1,2,4,5\}$ are VL-stable for every matrix in
$Q_A$, and the inequalities (K) hold, $r_1=1>0$, $r_2=9>3r_1$, $r_5<0$ and
$\lambda_5(F_{12})=\frac45>0$. The rate signs and $g>1$ are the conditions
(b) and (d) below. Every number is checked in exact arithmetic by
\texttt{hks\_certificates.jl}, which also finds $\rho(D^{-1}N)\approx1.10$.
\end{proof}

So the failure does not depend on large entries. We keep the matrix $A$ of
Theorem~\ref{thm:main} because its equilibria and rates are particularly
simple, which makes \S\ref{sec:proof} easy to check by hand.

We now describe the set of $r$ for which the proof of \S\ref{sec:proof}
applies verbatim, and estimate its size.

\paragraph{Which magnitudes the proof uses.} Let $\tilde A\in Q_A$.
\begin{itemize}
\item In each of the faces $\{1,3,4\}$, $\{1,2,3,4\}$, $\{1,2,3,5\}$ and
  $\{1,2,5\}$, the interactions are predator--prey pairs through species~1
  and at most one one-way link ($4\to2$ or $2\to5$). Hence $\tilde A_{KK}$ is
  VL-stable for every $\tilde A\in Q_A$. Choose $d_1=1$ and $d_k$ so that
  the predator--prey entries cancel; then the diagonal weight of the
  one-sided partner can be made large or small enough. The same block
  structure makes every feasible $F_k$ on $\Gamma$ internally stable.
\item In the face $K=\{1,2,4,5\}$, the argument of (d)(ii) uses only
  \[
    r_1>0,\qquad r_2>\frac{|\tilde a_{21}|}{|\tilde a_{11}|}\,r_1,\qquad
    r_5<0,\qquad \lambda_5(F_{12})>0.\tag{K}
  \]
  The bound $x_2\ge x_2(F_{12})$ at the relevant vertices holds for every
  $\tilde A\in Q_A$: adding species~4 raises $x_2$ ($\tilde a_{24}>0$), and
  removing the predator~1 raises it too.
\end{itemize}
So, for every $\tilde A\in Q_A$, the conclusion of
Theorem~\ref{thm:main} holds for every $r$ in the set $R(\tilde A)$ defined
by:
\begin{enumerate}[label=(\alph*)]
\item the interior equilibrium $x^*=-\tilde A^{-1}r$ is positive;
\item $F_{13},F_{134},F_{124},F_{125}$ are feasible, and their invasion
  rates have the signs of the table in (c);
\item the inequalities (K) hold, and all invasion rates of the face system on
  $K$ are nonzero;
\item the loop gain $g$ of (e) exceeds $1$.
\end{enumerate}
For fixed $\tilde A$, equilibria and invasion rates are linear in $r$.
Hence (a)--(c) define an open polyhedral cone with finitely many hyperplanes
removed, and (d) is the quartic inequality $\det C<0$. So $R(\tilde A)$ is
an open cone.

\paragraph{The cone for $\tilde A=A$.} Conditions (a)--(c) consist of 29
strict linear inequalities, besides the nonvanishing of the rates. Exact
enumeration shows that their cone is simplicial, with the five facets
\begin{gather*}
  r_1>0,\qquad r_5<0,\qquad 2r_1-3r_2+4r_3>0,\\
  2r_1+r_2-4r_3+12r_4>0,\qquad 72r_2-144r_1-r_4+24r_5>0,
\end{gather*}
namely $r_1>0$, $\lambda_5(F_{13})=r_5<0$, $\lambda_2(F_{13})<0$,
$\lambda_3(F_{124})<0$ and $\lambda_4(F_{125})<0$. Its extreme rays are
$(1,2,1,0,0)$, $(0,12,9,2,0)$, $(0,4,3,288,0)$, $(0,4,865,288,0)$ and
$(0,144,108,24,-431)$. In particular $x^*>0$ and the other two
inequalities of (K) follow from these five. So
\[
  R(A)=\{r:\ \text{the five inequalities hold and } g(r)>1\},
\]
minus the finitely many hyperplanes on which a rate of the face system on
$K$ vanishes. For example, every $r$ with $|r_i-(r_c)_i|\le\frac{68}{125}$
for all~$i$, where $r_c=(1,8,7,5,-7)$, satisfies the five inequalities and
$g(r)>1$; at $r_c$ the gain is $g\approx4.58$. This is checked exactly, with
interval bounds on~$g$. For the $r$ of Theorem~\ref{thm:main} the certified
radius is only $\frac{33}{20000}$, because $g$ is close to~$1$ there.

\paragraph{Its size.} We sampled $2\cdot10^5$ directions $r$ uniformly
on the unit sphere $S^4$ (script \texttt{hks\_rcone.jl}, which also tests
(K) on random examples in exact arithmetic).
\begin{itemize}
\item For $\tilde A=A$, $5.6\%$ of the directions give an interior
  equilibrium. $R(A)$ contains $0.08\%$ of all directions, which is
  $1.4\%$ of those with an interior equilibrium.
\item Within $R(A)$ the loop gain ranges up to about $10$. The $r$ of
  Theorem~\ref{thm:main} has $g\approx1.03$, so less marginal choices exist.
\item For $\tilde A$ with each magnitude of $A$ multiplied by a random
  factor in $[e^{-1},e]$, the proportion is $0.8\%$.
\end{itemize}

$R(\tilde A)$ is only the part where this particular cycle and this proof
apply. For a wider picture we ran the numerical census behind the remarks of
\S1 on $4000$ directions $r$ with an interior equilibrium, for $\tilde A=A$.
Its connections come from numerically integrated unstable manifolds; the
characteristic matrix is evaluated exactly.
\begin{itemize}
\item $8.2\%$ have an attracting simple heteroclinic cycle. There are four
  itineraries, and that of Theorem~\ref{thm:main} accounts for $3.8\%$.
\item $82.1\%$ have an acyclic invasion graph with nonzero rates. These
  systems are permanent \citep[Thm.~1]{hofbauer2022}, by the uniqueness of
  saturated equilibria for P-matrices.
\item In the remaining $9.7\%$ the only heteroclinic cycles found are
  repelling (or, in $0.1\%$, an integration did not converge), and this
  test does not decide permanence.
\end{itemize}
So, for this $A$, roughly one growth vector in twelve with an interior
equilibrium leads to an attracting heteroclinic cycle.

\paragraph{Positive growth rates.} All these growth vectors have $r_5<0$.
For $n\ge6$ this is not needed.

\begin{proposition}\label{prop:rpos}
Extend $A$ by a sixth species with $a_{66}=-1$ and $a_{56}=-3$ (species~6
harms~5), all other new entries zero, and let $r=(1,3,3,1,1,1)$. The
resulting $6\times6$ matrix satisfies (C1) and (C2), the system has the
interior equilibrium $(x^*,1)$, and it has an asymptotically stable
heteroclinic cycle through the equilibria $(F_I,1)$,
$I=\{1,3\},\{1,3,4\},\{1,2,4\},\{1,2,5\}$. Extending further as in
Corollary~\ref{cor:alln} gives counterexamples with $r\gg0$ for every
$n\ge6$.
\end{proposition}

\begin{proof}
The new link $6\dashv5$ is one-way, so no new nonzero cycle appears. At
every equilibrium containing species~6, $x_6=1$, and species~5 sees the
growth rate $1-3=-2$, so all equilibria and rates of
\S\ref{sec:proof} reappear. The faces used there, with~6 added, remain
VL-stable, since~6 is linked to them by at most one one-way link. In the face
$\{1,2,4,5,6\}$, species~6 is never lost: when absent its rate is $1>0$. The
argument of \S\ref{sec:proof}(d)(ii) then applies both to the vertices
containing~6 and to those without it, where species~5 has rate
$1+12x_2>0$ whenever absent. The script \texttt{hks\_certificates.jl}
checks every hypothesis of this proof in exact arithmetic; the face has 22
feasible equilibria, all with nonzero rates.
\end{proof}

The mechanism is that of Theorem~\ref{thm:main}. Species~6, present
throughout, lowers the effective growth rate of species~5 below zero, so the
loop still contains an effective obligate consumer. Positivity of $r$ is
therefore no obstruction for $n\ge6$, but it does not give a new mechanism.

For $n=5$ there is no counterexample with $r\gg0$
(Theorem~\ref{thm:five}(b), proved in the companion paper). The exact
sign-level tests there leave only two candidate itineraries with $r\gg0$,
among all 630 candidate patterns, and a single exact certificate excludes
both, for all magnitudes. So at $n=5$ every counterexample has an obligate
consumer.

\section{Discussion}\label{sec:disc}

\paragraph{Where the conjecture fails.} The conditions $a_{ii}<0$,
$a_{ij}a_{ji}\le0$ for $i\ne j$, and ``no nonzero cycles of length $\ge3$''
characterise qualitative VL-stability \citep[(1.2)]{hofbauer2008}. With an
interior equilibrium such systems are globally stable, hence permanent.
(C1) and (C2) allow, in addition, nonzero cycles of length $\ge3$ with
exactly one negative entry. Up to orientation, a nonzero 3-cycle of this
kind consists of two benefits and one harm. The counterexample uses one such
loop, $2\to5\dashv4\to2$, with all three effects one-sided.

The fifth species, a predator that can grow on its own, is present
throughout. It provides the background against which the loop switches its
members on and off (Figure~\ref{fig:orbit}). Each switch is an invasion. In
the characteristic matrix, the diagonal entries are the growth rates of the
invaders, and the negative entries are the decay rates of the other absent
species. The loop gain $g$ of \S\ref{sec:proof}(e) multiplies ratios of
decay to growth along the two cycles of dependencies. The heteroclinic cycle
attracts when decay wins, $g>1$.

\paragraph{Four species.} In a companion paper (in preparation) we prove, by
a computer-assisted argument in exact rational arithmetic, that
Conjecture~2.4 of \citet{hofbauer2008} holds for $n=4$: every $A$ with
$a_{ii}<0$ satisfying (C1) and (C2) is qualitatively permanent. The argument
has three parts.
\begin{itemize}
\item $-\tilde A$ is a P-matrix \citep[Lemma~6.2]{hofbauer2008}. So the
  system has a compact absorbing set and a unique saturated equilibrium
  \citep[Lemmas~6.3 and~6.4]{hofbauer2008}. An interior equilibrium is
  saturated, so it is that equilibrium; hence every boundary equilibrium has
  an invader.
\item For every $\tilde A\in Q_A$ and every $r$, the invasion graph has no
  cycle. Here zero invasion rates are allowed: an edge needs a positive rate
  for each gained species and a nonpositive rate for each lost one. The
  proof is an exact sign-level enumeration of candidate cycles, which are
  then excluded by certificates valid for all magnitudes.
\item \citet[Thm.~1]{hofbauer2022} then gives robust permanence. With zero
  rates, their nonvanishing assumption fails, and we extend their theorem:
  in a system with a P-matrix, a species with zero invasion rate cannot
  invade.
\end{itemize}
So $n=5$ is the smallest dimension in which the conjecture fails.

\paragraph{Five species: four motifs.} The same paper shows, by a
computer-assisted argument, that for $n=5$ a heteroclinic cycle needs at
least seven nonzero off-diagonal entries. Reducing all examples that we found to minimal patterns gives
four motifs. In the sign convention of $A$, they are
\begin{gather*}
  M_1=\begin{pmatrix}-&+&+&0&0\\-&-&0&+&0\\-&0&-&0&0\\0&0&0&-&-\\0&+&0&0&-\end{pmatrix},\qquad
  M_2=\begin{pmatrix}-&+&+&0&0\\0&-&+&+&0\\-&0&-&0&0\\0&0&0&-&-\\0&+&0&0&-\end{pmatrix},\\[4pt]
  M_3=\begin{pmatrix}-&+&+&0&0\\0&-&0&0&-\\-&0&-&-&0\\0&0&+&-&0\\+&0&0&0&-\end{pmatrix},\qquad
  M_4=\begin{pmatrix}-&0&+&0&0\\0&-&+&+&0\\-&0&-&+&+\\0&-&0&-&0\\0&0&-&0&-\end{pmatrix}.
\end{gather*}
$M_1$, $M_3$ and $M_4$ each contain exactly one nonzero cycle of length
$\ge3$: a 3-cycle of the kind above, with two further species attached by
predator--prey or one-way links. $M_2$ contains two such 3-cycles, which
share species~2. No motif has a nonzero cycle of length 4 or~5.
\begin{itemize}
\item $M_1$ is the pattern of $A$.
\item $M_2$ has the same loop. Species~1 preys on~3 and benefits one-sidedly
  from~2, and 2 benefits from~3. This closes a second loop,
  $3\to2\to1\dashv3$.
\item In $M_3$ the loop $1\to5\dashv2\to1$ contains the predator~1. Species
  3 is prey of~1, and 4 preys on~3.
\item $M_4$ has eight entries. In its loop $4\to3\to2\dashv4$ the harm is
  reciprocated, since 2 preys on~4. Species~3 is prey of~1 and preys on~5.
\end{itemize}

\begin{proposition}[computer-assisted]\label{prop:motifs}
For each $k=1,\dots,4$, the matrix $A^{(k)}$ and growth vector $r^{(k)}$
listed in \texttt{hks\_certificates.jl} have the sign pattern $M_k$, an
interior equilibrium and an asymptotically stable heteroclinic cycle.
Consequently, for every $n\ge5$, every $n\times n$ sign pattern with
negative diagonal that satisfies (C1) and (C2) and contains one of
$M_1,\dots,M_4$ is not qualitatively permanent.
\end{proposition}

\begin{proof}
For $M_1$ this is Theorem~\ref{thm:main}. For each example the script checks
in exact rational arithmetic, and prints, every hypothesis of the argument
of \S\ref{sec:proof}:
\begin{itemize}
\item (C1) and (C2), the P-matrix property, and $x^*>0$;
\item feasibility, internal stability (Routh--Hurwitz) and the signs of the
  invasion rates of every equilibrium on the cycle;
\item each connection, by Lemma~\ref{lem:conn} together with either
  Lemma~\ref{lem:VL} (with a diagonal VL certificate) or
  Lemma~\ref{lem:morse} (exact acyclicity of the face invasion graph, all
  rates nonzero, and a VL certificate for the target face);
\item simplicity, and a vector $q\gg0$ with $Cq\ll0$, i.e.\ $p=-q$ in
  Theorem~\ref{thm:hofbauer}.
\end{itemize}
Every hypothesis is a strict inequality, so the proofs of
Theorem~\ref{thm:robust} and Corollary~\ref{cor:pattern} apply.
\end{proof}

\paragraph{How the motifs were found.} The companion paper describes the
census in detail; here is an outline. We enumerated the $5\times5$ patterns
satisfying (C1) and (C2) up to relabelling. For each pattern we listed the
sign vectors of $r$ and the cyclic sequences of supports that could form a
heteroclinic cycle. Sign-level necessary conditions for edges of invasion
graphs, proved there for all magnitudes and all $r$ with the given signs,
exclude all but $630$ patterns ($4524$ sequences). This reduction is exact;
it replaces the continuum of $r$ by finitely many sign cases. For each
remaining sequence we searched magnitudes and $r$ for an attracting cycle,
maximising the spectral radius of \S\ref{sec:proof}(e), and we reduced every
example found by deleting interactions one at a time. The search reached an
attracting cycle in $598$ of the $630$ patterns, and each of them contains
one of $M_1,\dots,M_4$. This last step is a search; that nothing else occurs
is proved separately. In the $27$ patterns that contain none of the motifs,
each of the $113$ remaining sequences contains constraints that are
incompatible for all magnitudes and all $r$, by Gordan-type certificates in
exact arithmetic, some with case splits over sign conditions. Together with
the extension of \citet{hofbauer2022} to zero rates, this proves the
following.

\begin{theorem}[computer-assisted; proved in the companion paper]\label{thm:five}
Let $n=5$, and let $A$ have $a_{ii}<0$ and satisfy (C1) and (C2).
\begin{enumerate}[label=(\alph*)]
\item $A$ is qualitatively permanent if and only if its sign pattern
  contains none of $M_1,\dots,M_4$ after relabelling.
\item For every $\tilde A\in Q_A$ and every $r\gg0$ for which there is an
  interior equilibrium, the system is permanent.
\end{enumerate}
\end{theorem}

The ``only if'' direction of~(a) is Proposition~\ref{prop:motifs}, proved
here. Nothing in \S\S1--5 or in Proposition~\ref{prop:motifs} depends on
the companion paper or the census. What depends on them: the result for
$n=4$, so that $n=5$ is the smallest dimension; that seven entries are
needed; the rest of Theorem~\ref{thm:five}; and the census figures of
\S\ref{sec:rset}.

\paragraph{Biological implications.} We end with what these results mean
for ecological communities.

\emph{Signs are not enough.} Predicting the behaviour of a community from the
signs of its interactions alone has a long tradition: sign stability
\citep{quirk1965,may1973} and loop analysis \citep{levins1974} ask which
conclusions survive any choice of magnitudes. Qualitative permanence asks the
same about coexistence. (C1) restricts the pairs to consumer--resource
($+,-$), one-sided ($+,0$ or $-,0$) or absent interactions: no two species
compete directly or are mutualists. (C2) restricts the longer loops. For up
to four species these conditions guarantee coexistence for every choice of
magnitudes and of growth rates with an interior equilibrium (for $n\le3$ by
\citet{hofbauer2008}, for $n=4$ by our companion paper). From five species
on they do not. At $n=5$ the failure is localised: a sign pattern guarantees
coexistence exactly when it contains none of the four motifs $M_1,\dots,M_4$
(Theorem~\ref{thm:five}). Searching a sign digraph for these motifs is thus a
practical test. When a motif is present, the outcome depends on magnitudes
and growth rates; for the matrix of Theorem~\ref{thm:main}, about one growth
vector in twelve with an interior equilibrium gives an attracting
heteroclinic cycle (\S\ref{sec:rset}).

\emph{Invasibility is not coexistence.} In the counterexample, every
subcommunity that can persist on its own, i.e.\ every boundary equilibrium,
can be invaded by some missing species: the interior equilibrium is the only
saturated one \citep[Lemma~6.4]{hofbauer2008}. Criteria for coexistence
based on invasion growth rates \citep{chesson2000} therefore need more than
invasibility. \citet{hofbauer2022} show that invasibility suffices when the
invasion graph has no cycle, and \citet{spaak2023} discuss what this means
for coexistence theory. Theorem~\ref{thm:main} shows that this condition
cannot be dropped, even in the restricted class of (C1) and (C2): here the
invasions chain into a cycle (Figure~\ref{fig:graphs}). Species~4 joins
$\{1,3\}$, species~2 then displaces~3, species~5 displaces~4, and when~3
returns it displaces 2 and~5: after four invasions the community is back
where it started.
This is a cycle of community assembly in the sense of \citet{law1996}.

\emph{Long phases and deep troughs.} Along the heteroclinic cycle the
community spends ever longer periods close to one of four subcommunities,
$\{1,3\}$, $\{1,3,4\}$, $\{1,2,4\}$ and $\{1,2,5\}$, and the absent species
reach ever lower densities (Figure~\ref{fig:orbit}). An observer would see
long, apparently stable phases with some species missing, interrupted by
rapid turnover when a species re-invades from very low numbers. In the
deterministic model no species is ever lost. In a real population, densities
that correspond to a few individuals make extinction by demographic
stochasticity likely \citep{ovaskainen2010}, so after a few rounds one of the
loop species will in fact disappear. The heteroclinic cycle is then a
mechanism of species loss that a local stability analysis of the interior
equilibrium, or a check of invasibility, does not detect.

\emph{The ingredients.} The mechanism differs from the best-known
heteroclinic cycle in ecology, the competitive rock--paper--scissors cycle of
\citet{may1975}, which (C1) excludes. It needs three ingredients that are
common in natural communities.
\begin{itemize}
\item One-sided effects. Facilitation, where one species benefits from
  another at no cost to it, is widespread \citep{bruno2003}; amensalism, where
  one species harms another without being affected, occurs through
  interference or chemical inhibition. The loop $2\to5\dashv4\to2$ consists of
  two facilitations and one amensal effect. Reciprocal effects may be added
  as long as (C1) and (C2) hold (Corollary~\ref{cor:pattern}).
\item An obligate consumer. Species~5 cannot grow alone ($r_5<0$) and
  depends on species~2. At $n=5$ this is necessary: if every species can grow
  alone, the system is permanent (Theorem~\ref{thm:five}(b)). For $n\ge6$ the
  dependence can arise from the community instead, through a species that
  suppresses~5 (Proposition~\ref{prop:rpos}).
\item A generalist that persists throughout. Species~1 eats 2 and~3 and has
  resources of its own ($r_1>0$). It is never lost, and it sets the
  background against which the loop turns over.
\end{itemize}
As an illustration only, these roles fit a microbial community with a
grazer (species~1) that also has food of its own. Strain~5 cannot grow
alone and lives on a by-product of strain~2, which it does not affect.
Strain~5 also inhibits strain~4 with a toxin, and strain~4 improves the
conditions for strain~2.

\emph{Strong interactions.} The loop must be strong compared with
self-limitation: by Proposition~\ref{prop:loop}, a loop ratio below~$8$
makes the system globally stable. With the diagonal normalised to $-1$, the
three loop interactions in the example of Proposition~\ref{prop:normal}
have strength~$3$. This agrees with the view that strong self-regulation and
weak interactions stabilise communities \citep{may1973,mccann1998}. Our
results add that, below five species and in the class of (C1) and (C2), no
interaction strengths can destroy coexistence, and that from five species on
a single strong feedback loop of facilitation and harm can.

\paragraph{Open questions.}
\begin{itemize}
\item For $n\ge6$, does every pattern that satisfies (C1) and (C2) but is
  not qualitatively permanent contain one of $M_1,\dots,M_4$? Longer loops
  with one negative entry may give new motifs.
\item Is there a proof of the result for $n=4$, and of
  Theorem~\ref{thm:five}, without a computer?
\end{itemize}

\appendix
\section{Proof of Proposition~\ref{prop:boundary}}\label{app:boundary}

Let $\Gamma$ be simple, hyperbolic and IS, with equilibria $F_k$, supports
$I_k$, invaders $j_k$ and connections $\gamma_k$.

\emph{Structure.} As in the proof of Lemma~\ref{lem:conn}, $W^u(F_k)$ is a
curve. One branch of it is $\gamma_k$, which lies in $\relint X_{I_k\cup\{j_k\}}$;
the other branch leaves $X$. In particular $I_{k+1}\subseteq I_k\cup\{j_k\}$.

\begin{lemma}\label{lem:chain}
Let $G$ be a face of $X$, and $K\subset G$ a compact invariant set consisting
of finitely many equilibria $E$, hyperbolic within $G$, and finitely many
orbits connecting them. Assume:
\begin{itemize}
\item $W^u(e)\cap G\subseteq K$ for every $e\in E$;
\item the directed graph on $E$ given by these orbits has no cycle.
\end{itemize}
Then $K$ is asymptotically stable relative to $G$.
\end{lemma}

\begin{proof}
Induction on $|E|$.

\emph{Base case, $E=\{e\}$.} Then $W^u(e)\cap G=\{e\}$. Hartman--Grobman and
the $\lambda$-lemma \citep[Ch.~2]{palis1982} show that orbits starting in $G$
near $e$, but not on the local stable manifold, leave a fixed small ball
close to $W^u_{\mathrm{loc}}(e)$. Their exit points lie in $G$, which is
closed. If such points existed arbitrarily close to $e$, the exit points
would accumulate on a point of $W^u(e)\cap G\setminus\{e\}=\emptyset$, which
is impossible. So $e$ is asymptotically stable in $G$.

\emph{Induction step.} Some $e\in E$ has no incoming orbit. Removing $e$ and
the orbits leaving it gives a compact invariant set $K'$ satisfying the
hypotheses. By induction $K'$ is asymptotically stable in $G$. The compact
set of exit points from a small ball around $e$ along
$W^u_{\mathrm{loc}}(e)\cap G$ consists of points of $K$. These are carried
into a prescribed neighbourhood of $K'$ in finite time. By continuity of the
flow and the base-case argument, orbits from a small neighbourhood of $e$ in
$G$ either converge to $e$, or follow these orbits closely into the basin of
$K'$. This gives stability and attraction of $K$ in $G$.
\end{proof}

\begin{proof}[Proof of Proposition~\ref{prop:boundary}]
Fix a facet $G=\{x_j=0\}\cap X$. The set $K_G=\Gamma\cap G$ consists of the
$F_k$ with $j\notin I_k$, and of the $\gamma_k$ with
$j\notin I_k\cup\{j_k\}$.
\begin{itemize}
\item Each such $F_k$ is hyperbolic within $G$. Its unstable manifold within
  $G$ is either $\{F_k\}\cup\gamma_k$ (if $j_k\neq j$) or $\{F_k\}$ (if
  $j_k=j$), so it lies in $K_G$.
\item The graph of $K_G$ is a proper subgraph of the $m$-cycle, because
  species $j$ invades at some $F_p$ and then $\gamma_p\not\subset G$. So it
  is acyclic.
\end{itemize}
By Lemma~\ref{lem:chain}, $K_G$ is asymptotically stable in $G$.

A point of $\partial X$ near $\Gamma$ lies in some facet $G$, and it is close
to $K_G$, by compactness of $\Gamma$. Taking the minimum of the finitely many
neighbourhood sizes over the facets gives stability and attraction within
$\partial X$.
\end{proof}

\paragraph{Code availability.} Four Julia scripts (standard library only)
are included as ancillary files. \texttt{hks\_note\_checks.jl} redoes
every number of \S\ref{sec:proof} in exact rational arithmetic.
\texttt{hks\_rcone.jl} does the computations of \S\ref{sec:rset} on
$R(\tilde A)$ and the loop ratio. \texttt{hks\_certificates.jl} certifies
the heteroclinic cycles of Propositions~\ref{prop:rpos}
and~\ref{prop:motifs} in exact arithmetic, from the data in
\texttt{hks\_certificates\_\allowbreak data.jl}. \texttt{hks\_fill\_signs.jl}
lists the sign fillings of the zeros of $A$ that keep (C1) and (C2). All are
archived, with their expected output, at \zenodolink. The census code behind the remarks of \S1 and \S\ref{sec:disc}
and the census of \S\ref{sec:rset} will be published with the companion
paper, and is available from the authors.

\paragraph{Use of AI tools.} The authors used Claude (Anthropic), a large
language model, in this work. It was used to aid in the design and for
writing the computer searches and the exact-arithmetic verification scripts
(the ancillary files), to check the proposed proofs, and to draft and edit
the text. The authors checked all mathematical arguments and computations,
and take full responsibility for the content of this article.

\bibliographystyle{plainnat}
\bibliography{HKS_counterexample}

\end{document}